\documentstyle[12pt]{article}
\def\rbiprod{{\cdot\kern-.33em\triangleright\!\!\!<}}
\def\lbiprod{{>\!\!\!\triangleleft\kern-.33em\cdot\, }}

\def\lrbiprod{{\ \cdot\kern-.60em\triangleright\kern-.33em\triangleleft\kern-.33em\cdot\, }}

\def\lprod{{>\!\!\!\triangleleft\kern-.33em\ \, }}

\title{The Factorization of  Braided  Hopf Algebras  }
\author{Shouchuan Zhang, Yange Xu\\
 Department  of Mathematics, Hunan
University\\ Changsha  410082, \
 P.R. China.
  }
\date{}
\begin{document}
\newtheorem{Theorem}{\quad Theorem}[section]
\newtheorem{Proposition}[Theorem]{\quad Proposition}
\newtheorem{Definition}[Theorem]{\quad Definition}
\newtheorem{Corollary}[Theorem]{\quad Corollary}
\newtheorem{Lemma}[Theorem]{\quad Lemma}
\newtheorem{Example}[Theorem]{\quad Example}
\newtheorem{Remark}[Theorem]{\quad Remark}

\maketitle \addtocounter{section}{0}

\begin {abstract}
We obtain   the double factorization  of braided bialgebras or
braided Hopf algebras,  give relation among integrals and
semisimplicity of
 braided Hopf algebra  and    its factors.

\end {abstract}

\section*{\bf Introduction}

It is well-known  that  the factorization of domain  plays an
important role in ring theory. S. Majid in \cite [Theorem 7.2.3]
{Ma95b} studied the factorization  of Hopf algebra and showed that
$H \cong A \bowtie B$ for two sub-bialgebras $A$ and $B$ when
multiplication $m_H$ is bijective. C. S. Zhang, B.Z. Yang  and B. S.
Ren \cite {ZYR05} generalized these results to braided cases.

Braided tensor categories become more and more important. They have
been  applied in conformal field, vertex operator algebras, isotopy
invariants of links (See \cite {Hu05, HK04, BK01}, \cite {He91,Ka97,
Ra94}).

In this paper, we  obtain  the double factorization  of braided
bialgebras or braided Hopf algebras, i.e. we give the conditions to
factorize a braided Hopf algebra into the double cross products of
sub-bialgebras or sub-Hopf algebras. We give relation among
integrals and semisimplicity of
 braided Hopf algebra  and    its factors.

Throughout this paper, we work in braided tensor category $({\cal
C}, C)$, where ${\cal C}$ is a concrete category and  underlying set
of every  object in  ${\cal C}$ is a  vector space over a field $k$.
For example, Yetter-Drinfeld category over Hopf algebras with
invertible antipode and some important categories in \cite {BK01}
are such categories.

\section*{\bf Preliminaries}\label{s0}

We assume  that   $H$ and $A$ are
 two braided bialgebras  with  morphisms:
\begin {eqnarray*}
\alpha : H \otimes A \rightarrow &A& , \hbox { \ \ \ \ }
\beta : H \otimes A \rightarrow H,    \\
\phi :  A \rightarrow H \otimes  &A&   , \hbox { \ \ \ \ } \psi : H
\rightarrow H \otimes A
\end {eqnarray*}
                        such that $(A, \alpha )$ is a left $H$-module coalgebra,
 $(H, \beta )$ is a right $A$-module coalgebra,
 $(A, \phi )$ is a left $H$-comodule algebra, and
 $(H, \psi )$ is a right $A$-comodule algebra.

We define   the multiplication $m_D$ , unit $\eta _D$,
comultiplication  $\Delta _D$ and counit $\epsilon _D$  in $A
\otimes H$  as follows:
$$\Delta _D = (id _A \otimes m_H \otimes m_A \otimes id _H )(id _A \otimes id _H \otimes C_{A, H} \otimes id _A \otimes id _H )
(id _A \otimes \phi  \otimes \psi \otimes id _H )(\Delta _A \otimes
\Delta _H ), \ \ \ $$
$$m _D =(m _A \otimes
m _H )(id _A \otimes \alpha   \otimes \beta \otimes id _H )(id _A
\otimes id _H \otimes C_{H, A} \otimes id _A \otimes id _H ) (id _A
\otimes \Delta _H \otimes \Delta_A \otimes id _H ) , \ \ \ $$ and
$\epsilon _D = \epsilon _A \otimes \epsilon _H$ , $ \eta _D = \eta
_A \otimes \eta _H.$ We denote   $(A \otimes H, m_D, \eta _D, \Delta
_D, \epsilon _D) $
  by                      $ A ^{\phi} _{\alpha}  {\bowtie }
^{\psi}_{\beta} H ,$  which is
 called the
double bicrossproduct of $A$ and $H$. We  denote it by  $A \stackrel
{b} {\bowtie }  H$ or $A\lrbiprod  H$ in short.
 When $\phi $  and $\psi$ are
trivial,  we   denote $ A_\alpha ^\phi \bowtie _\beta ^\psi H $  by
 $A _\alpha {\bowtie } _\beta H$ or  $A  {\bowtie }  H$, called  a double cross product
 (see \cite {ZC01, BD99, Zh99}).

 ${\cal V}ect(k)$  denotes the braided tensor category of all vector
spaces over  field $k$, equipped with ordinary tensor and unit
$I=k$, as well as ordinary  twist  map as braiding. $^B_B {\mathcal
YD}$ denotes Yetter-Drinfeld category (see \cite [Preliminaries]
{Zh99}).

\section {The factorization of  braided  bialgebras or braided Hopf algebras
} \label {s1}

In this section, we obtain the factorization  of  braided bialgebras
or braided Hopf algebras.

The associative law   does not hold for double cross products in
general, i.e. equation
$$  (({ A_1 } { \ }  _{\alpha _1 } \bowtie _{\beta _1} A_2 )
{ \ }_{\alpha _2}  \bowtie _{\beta _2}  A_3) =
  ( {A_1 } { \ } _{\alpha _1 } \bowtie _{\beta _1} ({A_2}
  { \ } _{\alpha _2}  \bowtie _{\beta _2}  A_3))$$
  does not holds in general. Therefore we denote a method of adding brackets for $n$ factors by $\sigma $.
For example, when $n=3$, $$\sigma_{1}( { A_1 } { \ }  _{\alpha _1 }
\bowtie _{\beta _1} A_2 { \ }_{\alpha _2}  \bowtie _{\beta _2}
A_3)=(({ A_1 } { \ }  _{\alpha _1 } \bowtie _{\beta _1} A_2 ) { \
}_{\alpha _2}  \bowtie _{\beta _2}  A_3)$$ and $$\sigma_{2}( { A_1 }
{ \ }  _{\alpha _1 } \bowtie _{\beta _1} A_2 { \ }_{\alpha _2}
\bowtie _{\beta _2} A_3)= ( {A_1 } { \ } _{\alpha _1 } \bowtie
_{\beta _1} ({A_2}
  { \ } _{\alpha _2}  \bowtie _{\beta _2}  A_3))$$

\begin {Definition}  \label {4.2.1}
   Let $A_1, A_2, \cdots , A_n$ and $H$ be braided bialgebras or braided Hopf algebras,
   $\sigma $ be a method adding brackets for n factors. If $j_{A_1}
(A_1) j_{A_2} (A_2) \cdots j_{A_n} (A_n)= H$ and
$j_{A_i}:A_{i}\rightarrow H$ is a braided bialgebra or Hopf algebra
morphism for $i = 1, 2, \cdots , n.$,  If  for every pair of
brackets of $\sigma $ :
$$(A_{l} \otimes A_{l+1} \otimes \cdots \otimes A_{l+t}),$$
 $j_{A_l}(A_l)  j_{A_{l+1}} (A_{l+1})\cdots  j_{A_{l+t}}(A_{l+t})$ is
 a sub-bialgebra or sub-Hopf algebra of $H$ \ \  \ and  \\
 $m_H ^{t} (j_{A_l} \otimes j_{A_{l+1}} \otimes
\cdots \otimes j_{A_{l+t}})$  is a bijective map from $A_{l} \otimes
A_{l+1} \otimes A_{l+2} \otimes \cdots \otimes A_{l+t}$  \ \ \ \ \
onto  \\
$j_{A_l}(A_l)  j_{A_{l+1}} (A_{l+1})\cdots  j_{A_{l+t}}(A_{l+t})$,
then $\{ j_{A_1}, j_{A_2}, \cdots , j_{A_n} \}$    is called a
double factorization of $H$ with respect to $\sigma$. If for every
method $\sigma$  adding brackets,
 $\{ j_{A_1}, j_{A_2}, \cdots , j_{A_n} \}$
a double factorization of $H$ with respect to $\sigma$, then
 $\{ j_{A_1}, j_{A_2}, \cdots , j_{A_n} \}$    is called
a double factorization of $H$. If $A_t$ is a sub-bialgebra or
sub-Hopf algebra of $H$ and $j_{A_t}$  is a inclusion map from $A_t
$ to $H$ by sending $a$ to $a$ for any $a\in A_t$,  $t =1, 2, \cdots
, n$, then $H = A_1A_2 \cdots A_n$ is called an inner double
factorization of $H$.
\end {Definition}

\begin {Theorem}  \label {4.2.2}
   Let $A_1, A_2, \cdots , A_n$ and $H$ be braided bialgebras
or braided Hopf algebras, and let $\sigma $ be a method adding
brackets for n factors. Assume that $j_{A_i}$  is  a bialgebra or
Hopf algebra morphism from $A_i$ to $H$ for $i = 1, 2, \cdots , n.$
If
 $\{ j_{A_1}, j_{A_2}, \cdots , j_{A_n} \}$
a double factorization of  $H$ with respect to  $\sigma$
 then,                 in braided tensor category
$({\cal C}, C)$, there exists   a set  $\{ \alpha _i, \beta _i \mid
i= 1, 2, \cdots , n\}$ of morphisms    such that
 $$ \sigma ( {A_1} { \ } _{\alpha _1 } \bowtie _{\beta _1} {A_2 } { \ }
_{\alpha _2} \bowtie _{\beta _2} {A_3} { \ }_{\alpha _3} \bowtie
_{\beta _3} \cdots _{\alpha _{n-1}}\bowtie  _{\beta _{n-1}}A_n )
\cong H
 \hbox { \ \ \   (   as bialgebras or Hopf algebras ) }$$  and
 the isomorphism is $m^{n-1}_H (j_{A_1} \otimes j_{A_2} \otimes
 \cdots \otimes j_{A_n})$.

\end {Theorem}
{\bf Proof.}  We use induction for $n$.  When $n=2$, we can obtain
the proof by \cite [Theorem 2.1]{ZYR05} (i.e. the factorization
theorem). For $n> 2$, we can assume  that
$$ \sigma ( A_1  \otimes  A_2
\otimes A_3 \cdots  \otimes A_n ) = \sigma _1 ( A_1 \otimes  A_2
\otimes \cdots \otimes
  A_t )
\otimes  \sigma _2 ( A_{t+1}  \otimes  \cdots \otimes A_n ). $$

Next we consider $t$ in following three cases.

 (i)\ \ If $1<t<n-1$, let $B_1 = j_{A_1} (A_1)j_{A_2}(A_2)\cdots
j_{A_t}(A_t)$ \ \ \ and \\
$B_2 = j_{A_t+1}(A_{t+1})j_{A_{t+2}} (A_{t+2}) \cdots j_{A_n}
(A_n).$  It follows from the inductive assumption that  $\{ j_{A_1},
j_{A_2}, \cdots , j_{A_t} \}$ is a double factorization of $B_1$
with respect to $\sigma _1$, $\{ j_{A_{t+1}}, j_{A_{t+2}}, \cdots ,
j_{A_n} \}$ is a double factorization of $B_2$ with respect to
$\sigma _2$. Thus, there exists a set  $\{\alpha _i, \beta _i \mid
 i= 1, 2, \cdots  ,n, i\not= t  \} $ of morphisms such that
\begin {eqnarray}\label {le13e1} \sigma _1 ( {A_1} { \ } _{\alpha _1 } \bowtie _{\beta _1} {A_2 }
{ \ } _{\alpha _2} \bowtie _{\beta _2}  \cdots {\alpha _{t-1}}
\bowtie _{\beta _{t-1}} A_t ) \cong B_1 \hbox {  (   as bialgebras
or Hopf algebras ) }
 \end {eqnarray}
 and
\begin {eqnarray}\label {le13e2} \sigma _2( {A_{t+1}} { \ }_{\alpha _{t+1}} \bowtie _{\beta _{t+1}} \cdots
_{\alpha _n}\bowtie  _{\beta _n}A_n ) \cong B_2 \ \hbox {  (   as
bialgebras or Hopf algebras ) }.\end {eqnarray}
 The isomorphisms of (\ref {le13e1}) and (\ref {le13e2})  are
 $m^{t-1}_H (j_{A_1} \otimes j_{A_2} \otimes
 \cdots \otimes j_{A_t})$ and
  $m^{n-t}_H (j_{A_{t+1}} \otimes j_{A_{t+2}} \otimes
 \cdots \otimes j_{A_n})$, respectively.
Let $j_{B_1}  $  and $j_{B_2}$ denote the contain-map  of $B_1$  and
$B_2$ in $H$ respectively.  we can get that $m_H (j_{B_1} \otimes
j_{B_2})$ is a bijective since
 $m^{n-1}_H (j_{A_1} \otimes j_{A_2} \otimes
 \cdots \otimes j_{A_n})$ is a bijective.
In fact, $m_H (j_{B_1} \otimes j_{B_2}) = m^{n-1}_H (j_{A_1} \otimes
j_{A_2} \otimes
 \cdots \otimes j_{A_n}) ((m^{t-1}_H (j_{A_1} \otimes j_{A_2} \otimes
 \cdots \otimes j_{A_t}))^{-1} \otimes (m^{n-t-1}_H (j_{A_{t+1}} \otimes
 j_{A_{t+2}} \otimes
 \cdots \otimes j_{A_n}))^{-1}).$

(ii)\ \ If $t=1$, let $B_{1}=A_{1}$, $B_{2}=j_{A_2}(A_{2})j_{A_{3}}
(A_{3}) \cdots j_{A_n} (A_n)$,  let $j_{B_{1}}=j_{A_{1}}$ and let
$j_{B_{2}}$ be  a contain-map from $B_{2}$ to $H$. We can  get that
$m_H (j_{B_1} \otimes j_{B_2})$ is a bijective by the way similar to
(i).

(iii)\ \ If $t=n-1$, let $B_{1}=j_{A_1}(A_{1})j_{A_{2}} (A_{2})
\cdots j_{A_{n-1}} (A_{n-1})$ and  $B_{2}=A_{n}$. We can also  get
that $m_H (j_{B_1} \otimes j_{B_2})$ is a bijective.

Consequently,   by \cite [ Factorization Theorem ]{ZYR05},
 there exist $\alpha _t$  and $\beta _t$ such that $ {B_1} \
_{\alpha _t } \!\bowtie _{\beta _t} B_2  \cong
j_{B_{1}}(B_1)j_{B_{2}}(B_2)=H.$ Considering relations (\ref
{le13e1}) and (\ref {le13e2}), we have that
  $$ \sigma ({ A_1}{ \ }  _{\sigma _1 } \bowtie _{\beta _1} {A_2} { \ }
_{\alpha _2} \bowtie _{\beta _2}{ A_3} { \ } _{\alpha _3} \bowtie
_{\beta _3} \cdots _{\alpha _{n-1}}\bowtie  _{\beta _{n-1}}A_n )
\cong H
 \hbox { \ \ \   (   as bialgebras or Hopf algebras ) }$$
 and the isomorphism is
 $m^{n-1}_H (j_{A_1} \otimes j_{A_2} \otimes
 \cdots \otimes j_{A_n})$.
We complete the proof. $\Box$

\begin {Corollary}  \label {4.2.3}
   Let $X$, $A$ and $H$ be braided bialgebras
or braided Hopf algebras. Assume $j_{A}$ and $j_H$  are  bialgebra
or Hopf algebra morphisms
 from $A$ to $X$  and $H$ to $X$, respectively.
 Then
 $\{ j_{A}, j_{H} \}$  is
a double factorization of $X$ iff                  in braided tensor
category $({\cal C}, C),$ there exist morphisms $ \alpha  $ and
$\beta $
    such that
 $$ {A} { \ } _{\alpha } \bowtie _{\beta } {H }
 \cong X
 \hbox { \ \ \   (   as bialgebras or Hopf algebras ) }$$  and
 the isomorphism is $m_X (j_{A} \otimes j_{H}).$

\end {Corollary}

\begin {Corollary}  \label {4.2.4}
   Let $A_1, A_2, \cdots , A_n$ be braided sub-Hopf algebras of a finite-dimensional Hopf algebra
   $H$, and $H=A_1 A_2 \cdots  A_n$.
Assume that  $\sigma $  is a method adding brackets for $n$ factors.
Then the following statements are equivalent.

(i) $H= A_1A_2 \cdots A_n$ is an inner double factorization of $H$
with respect to $\sigma $.

(ii)   $dim  H = dim  (A_1) dim  (A_2) \cdots dim  (A_n)$, and for
every pair of brackets in $\sigma $ :   $(A_t \otimes A_{t+1}
\otimes \cdots \otimes A_{t+l} ),  $ $(A_t A_{t+1} \cdots  A_{t+l} )
$  is a sub-Hopf algebra of $H$.

(iii) For every pair of bracket in $\sigma $ :   $(A_t \otimes
A_{t+1} \otimes \cdots \otimes A_{t+l} ),  $ $(A_t A_{t+1} \cdots
A_{t+l} )  $  is a sub-Hopf algebra of $H$, and there exists   a set
$\{ \alpha _i, \beta _i \mid i= 1, 2, \cdots , n\}$  of morphisms
such that
 $$ \sigma ( {A_1} { \ } _{\alpha _1 } \bowtie _{\beta _1} {A_2 } { \ }
_{\alpha _2} \bowtie _{\beta _2} {A_3} { \ }_{\alpha _3} \bowtie
_{\beta _3} \cdots _{\alpha _{n-1}}\bowtie  _{\beta _{n-1}}A_n )
\cong H
 \hbox { \ \ \   (   as bialgebras or Hopf algebras ) }$$  and
 the isomorphism is $m^{n-1}_H $.


\end {Corollary}

{\bf Proof.} Obviously,  (iii) implies (ii).
 By Theorem \ref {4.2.2},
 (i) implies  (iii). It is sufficient to show that  (ii)
implies  (i). Assume that $(A_t \otimes A_{t+1} \otimes \cdots
\otimes A_{t+l})$ is  a pair of brackets in $\sigma $. We only need
to show that $m_H ^{l}$ is a bijective map from $(A_t \otimes
A_{t+1} \otimes \cdots \otimes A_{t+l})$ onto $A_t A_{t+1}\cdots
A_{t+l}$. Since

$ dim  ( A_1A_2\cdots A_n )= \ dim  (A_1) dim  (A_2) \cdots
        dim  (A_n),$
        we have that    $ dim  (A_t \otimes A_{t+1} \otimes \cdots \otimes
A_{t+l}) = dim  (A_t A_{t+1}\cdots A_{t+l})$, which implies $m_H
^{l}$ is a bijective map from $(A_t \otimes A_{t+1} \otimes \cdots
\otimes A_{t+l})$ onto $A_t A_{t+1}\cdots A_{t+l}$. We complete the
proof. $\Box$

\begin {Lemma}  \label {4.2.5}

Let $A$  and $B$ be braided sub-Hopf algebras of
 braided Hopf algebra $H$.

(i) If $AB=BA$,   then  $AB$   is sub-Hopf algebra of $H$.

(ii) If the antipode of $H$ is  invertible, and  $AB$ or $BA$ is a
braided sub-Hopf algebra of $H$, then $AB= BA$.
\end{Lemma}

{\bf Proof.} (i) It is clear.

(ii) We can assume that $AB$ is a braided sub-Hopf algebra of $H$
without lost generality.  For any $a \in A, b \in B, $  we see that
$S(S^{-1}(a) S^{-1}(b)) = ba$. Thus   $BA \subseteq AB$ since $AB$
is a braided sub-Hopf algebra.    For any $x \in AB$, there exist
$a_i \in A, b_i \in B$ such that
 $S(x) = \sum a_ib_i$. We see that
 $x = S^{-1} S(x) =
 S^{-1}(\sum a_ib_i) =
 \sum S^{-1} (b_i) S^{-1} (a_i) \in BA$. Thus $AB \subseteq BA$.
Consequently, $AB=BA$. $\Box$

\begin {Corollary}  \label {4.2.6}
   Let $A_1, A_2, \cdots , A_n$ be braided sub-Hopf algebras of  a finite-dimensional braided  Hopf algebra $H$, and $H=A_1 A_2 \cdots  A_n$.
   Then the following statements are equivalent.

(i) $H= A_1A_2 \cdots A_n$ is an inner double factorization of $H$.

(ii)   $dim  H = dim  (A_1) dim  (A_2) \cdots dim  (A_n)$ and
$A_uA_v = A_v A_u$  \ \ for  $1 \le u < v \le n$.

  (iii)  For $1\le u < v \le n$, $A_uA_v=A_vA_u$,
  and  for any method $\sigma $  adding brackets,
  there exists   a set  $\{ \alpha _i, \beta _i
\mid i= 1, 2, \cdots , n\}$  of morphisms   such that
 $$ \sigma ( {A_1} { \ } _{\alpha _1 } \bowtie _{\beta _1} {A_2 } { \ }
_{\alpha _2} \bowtie _{\beta _2} {A_3} { \ }_{\alpha _3} \bowtie
_{\beta _3} \cdots _{\alpha _{n-1}}\bowtie  _{\beta _{n-1}}A_n )
\cong H
 \hbox { \ \ \   (   as  Hopf algebras ) }$$  and
 the isomorphism is $m^{n-1}_H $.

(iv) $H=A_{i_1}A_{i_2} \cdots A_{i_n}$ is an inner double
factorization of $H,$ where $ \{ i_1, i_2, \cdots , i_n \} = \{ 1,
2, 3, \cdots , n\} $ as set.

\end {Corollary}

{\bf Proof.} $(i) \Rightarrow (ii)$  For  $ 1\le u < v \le n$, there
exists  a method adding brackets $\sigma $ such that $(A_u \otimes
A_v)$   is a pair of brackets in $\sigma $. Thus $A_uA_v$ is a
sub-Hopf algebra of $H$. By Lemma \ref {4.2.5}, $A_uA_v = A_vA_u.$
It follows from Corollary \ref {4.2.4} that
  $dim  H = dim  (A_1) dim  (A_2) \cdots dim  (A_n)$.

$(ii) \Rightarrow (i)$   follows from Corollary \ref {4.2.4}.

Similarly, (ii) and (iv)  are equivalent.

By Corollary \ref {4.2.4}, we also have that $(ii)$  and (iii)  are
equivalent. $\Box$

If $H$ is an almost commutative braided Hopf algebra, in particular,
$H$ is a coquasitriangular braided  Hopf algebra, then
 $AB = BA$  for any  braided sub-Hopf algebras $A$  and $B$ of $H$.
 Note that every quantum commutative braided Hopf algebra $H$ is a
coquasitriangular  braided Hopf algebra with coquasitriangular
structure $r = \epsilon _H \otimes \epsilon _H.$

\begin {Example} \label {4.9} (\cite[Lemma 3.4]{AS98}) Assume that $\Gamma$ is a commutative group and
$\hat \Gamma$ is the character group of $\Gamma$ with $g_i \in
\Gamma,$ $\chi _i\in \hat \Gamma$, $\chi _i (g_j)\chi_j (g_i) =1$,
$1< N_i $, where $N_i$ denotes the order of $\chi _i (g_i)$,  \  $1
\le i <  j \le \theta$. Let $H$ denote the algebra generated by set
$ \{ x_i \mid 1 \le i \le \theta \}$ with relation:
\begin {eqnarray} \label {qlse1} x_l ^{N_l}=0,\  x_{i}x_{j} = \chi_{j}(g_{i})
x_{j}x_{i}  \ \ \ \hbox { for }  1\le i, j, l \le \theta \hbox {
with } i \not= j.
\end {eqnarray} Define coalgebra operations and $kG$-(co-)module operations in
 $H$ as follows:
$$\Delta x_i = x_i \otimes 1 + 1 \otimes x_i, \ \ \epsilon (x_i)
=0,$$
 $$\delta ^-(x_{i}) = g_{i}\otimes x_{i}, \qquad h \cdot x_{i} =
\chi_{i}(h)x_{i}. $$ Then  $H$ is  called a quantum linear space in
$^{k\Gamma}_{k\Gamma} {\cal YD} $. By \cite [Lemma  3.4] {AS98}, $H$
is a braided Hopf algebra with  $dim H =N_{1}N_{2}\cdots
N_{\theta}$.
 Let
$H_{i}$ is the  sub-algebra  generated by $x_{i}$ in $H$. It is easy
to check that $H_i$ is a braided sub-Hopf algebra of $H$ with $dim
H_{i}=N_i$ and $H = H_1 H_2 \cdots H_\theta $ is an inner double
factorization of $H$ by Corollary \ref {4.2.6} (ii). Furthermore,
when $\theta =1$ and $N_1 =p$ is prime, then  $H$ is not commutative
with $dim H = p$.
\end {Example}

By the way, it is well-known that the 8th Kaplansky's conjecture is
that if the dimension of Hopf algebra $H$ is prime then $H$ is
commutative and cocommutative. Y. Zhu \cite {Zh94} gave the positive
answer. Now the  example above  show that braided version of the 8th
Kaplansky's conjectre  does not hold, i.e. there exists a
noncommutative  braided  Hopf algebra $H$ with prime dimension.

  If there are two  non-trivial  sub-bialgebras or sub-Hopf algebras
   $A$ and $B$ of $H$
 such that $H=AB$ is an inner double factorization of $H$,
then  $H$ is called a double factorizable bialgebra or Hopf algebra,
 Otherwise,
 $H$ is called a double infactorisable bialgebra or Hopf algebra.

A bialgebra (Hopf algebra) $H$  is said to satisfy the a.c.c. on
 sub-bialgebras (sub-Hopf algebras) if for every chain $A_1 \subseteq A_2
 \subseteq  \cdots $  of sub-bialgebras (sub-Hopf algebras) of $H$ there is an integer $n$
  such that $A_i = A_n$, for all $i > n.$  Similarly, we can define d.c.c..

   If $H$
 satisfies the d.c.c. or a.c.c. on sub-bialgebras (sub-Hopf algebras), then
  $H$ can be factorized  into a product of finite double infactorisable
 sub-bialgebras (sub-Hopf algebras).

  By the Corollary \ref {4.2.3},
  we can easily know that Sweedler's 4-dimensional Hopf algebra $H_4$
 over field $k$ is double infactorisable  in category ${\cal V}ect (k)$. In fact, if $H_4$ is double factorisable, then
  there exist two non-trivial sub- Hopf algebras $A$ and $B$ such that
  $H=AB$ and $ H \cong A \bowtie B$. It is clear that $A$ and $B$ are 2-dimensional. Thus they are
  commutative, which implies $H_4$  is commutative. We have a contradiction.
Thus $H_4$ is double infactorisable. Similarly, if $p$ and  $q$ are
two prime numbers and $H$ is a non-commutative Hopf algebra with
$dim  H = pq$, then $H$ is double infactorisable. Consequently,
every  Taft algebra $H$ with $dim  H= p^2$ is double infactorisable
in ${\cal V}ect (k)$.

\section {The factorization of  braided  bialgebras or braided Hopf algebras in Yetter-Drinfeld categories }\label {s1}
Throughout  this section, we work in  a Yetter-Drinfeld  category
$^B_B {\cal YD}$. In this section, we give relation among integrals
and semisimplicity of
 braided Hopf algebra  and    its factors.

If $ H$ is a finite-dimensional braided Hopf algebra in
$_{B}^{B}\mathcal{YD}$, then $\int _H^l$  and $\int _H^r$ are
one-dimensional by \cite
 {Ta99} or \cite [Theorem 2.2.1]{Zh99}, so there exist a non-zero left integral $\Lambda _H^l$ and  a non-zero
right integral $\Lambda _H^r$.

\begin {Proposition} \label {4.3.3} If $A, H$ and $D=A \stackrel {b} {  \bowtie } H$ are finite dimensional braided  Hopf
algebras, then

(i)  there are $u\in H, v\in A$ such that
 $$\Lambda _D {}^l =
 \Lambda _A^l \otimes u, \hbox { \ \ \ }
\Lambda _D^r = v \otimes \Lambda _H^r ; $$

(ii) If  $ D=A \stackrel {b} {  \bowtie } H$ is semisimple, then $A$
and $H$ are semisimple;

(iii) If  $ D=A \stackrel {b} {  \bowtie } H$ is unimodular, i.e.
$\Lambda_ D ^l = \Lambda_ D^r$,  then there exists a non-zero $x \in
k$ such that
             $ \Lambda_ D =x
 \Lambda _A ^l \otimes \Lambda _H^r,$
 and              $A \stackrel {b} {  \bowtie } H$  is semisimple iff
 $A$ and $H$ are semisimple.
 \end {Proposition}
 {\bf Proof .}  (i)
 Let $a  ^{(1)},  a  ^{(2)}, \cdots, a  ^{(n)}$  and $ h ^{(1)},  h  ^{(2)},
 \cdots ,  h^{(m)}$ be the basis of $A$  and  $H$,  respectively.
Assuming $$ \Lambda _D^l = \sum k_{ij} (a ^{(i)} \otimes h^{(j)}) $$
where $k_{ij} \in k$, we have that

$$a \Lambda _D^l = \epsilon (a) \Lambda _D^l $$ and
$$   \sum \epsilon (a) k_{ij} (a ^{(i)} \otimes h^{(j)})=
 \sum k_{ij} (a a ^{(i)} \otimes h^{(j)}) ,$$   for any $a \in A.$
Let $x_j = \sum _i k_{ij}a  ^{(i)}$. Considering $\{h^{(j)} \}$ is a
base of $H$,
 we get $ x_j $ is a left integral of $A$  and there exists  $k_j \in k$ such that
 $x_j = k_j \Lambda _A^l$  for $j = 1, 2, \cdots m$.
Thus $$\Lambda _D^l = \sum _j k_j (\Lambda  _A^l \otimes h^{(j)}) =
\Lambda _A^l \otimes  u, $$

where $u = \sum _j k_j h^{(j)}.$

Similarly, we have  that $\Lambda _D^r = v \otimes \Lambda _H^r$.

(ii) and (iii)  follow from  (i). $\Box$

\begin{Remark} For braided version of  biproduct $A\lbiprod H$ (
see \cite [Corollary 2.17] {BD99} or \cite [Chapter 4] {Zh99}),
Proposition \ref {4.3.3} also holds. For example, in Example \ref
{4.9}, $\int_H ^l = \int_H^r = k x_1 ^{N_1-1}x_2 ^{N_2-1}\cdots
x_\theta ^{N_\theta -1} $ of quantum linear space $H$. By
Proposition \ref {4.3.3}, the
 integral of the biproduct $D = H \lbiprod k\Gamma $ is
 $\int _D = k x_1 ^{N_1-1}x_2 ^{N_2-1}\cdots x_\theta  ^{N_\theta -1} \otimes
 (\sum _{g \in \Gamma}g)$.
\end {Remark}

\section {The factorization  of   ordinary Hopf algebras }\label {s1}

Throughout  this section, we work in  the category ${\cal V}ect (k)
$. In this section, using the results in preceding section, We give
relation among semisimplicity of
 Hopf algebra  and    its factors.

\begin {Lemma}  \label {4.3.4'} Assume that  $H$ is a finite-dimensional Hopf algebra.
Then

(i) $H$ is cosemisimple iff $H^*$ semisimple.

(ii)  $H$ is semisimple iff $H^*$ cosemisimple.

\end {Lemma} {\bf Proof.} (i) If $H$ is cosemisimple, then there exists $T \in \int
_{H^*}^l $ such that $\epsilon _{H^*} (T) \not= 0$ by \cite [Theorem
2.4.6, i.e. dual Maschke theorem] {Mo93}. Therefore $H^*$ is
semisimple. Conversely, if $H^*$ is semisimple, then there exists $T
\in \int _{H^*}^l $ such that $\epsilon _{H^*} (T) \not= 0$. Using
again \cite [dual Maschke theorem] {Mo93}, we have that $H$ is
cosemisimple.

Similarly,  we get (ii). $\Box$

\begin {Lemma} \label {4.3.1}  (cf.  \cite  {Mo93}, \cite  {Ra94} ).
If $H$ is a finite dimensional Hopf algebra  with $char  k =0 ,$
then the following conditions are equivalent.

(i) $H$ is semisimple.

(ii) $H$ is cosemisimple.

(iii) $H$ is semisimple and cosemisimple.

(iv) $S_H^2 = id _H$.

(v) $tr(S_H^2) \not=0.$
\end {Lemma}
{\bf Proof.}  $(i) \Rightarrow (ii)$. If  $H$ is semisimple, then
$H^*$ is cosemisimple by Lemma \ref {4.3.4'}, so  $H^*$ is
semisimple by \cite [Theorem 2.5.2] {Mo93}. Therefore $H$ is
cosemisimple.  It  follows from \cite [Theorem 2.5.2] {Mo93} that
(ii) implies (i). Consequently, (i), (ii) and (iii) are equivalent.
Using formula $tr (S_H^2) = \epsilon _H (\Lambda _H^l) \Lambda
_{H^*} ^r (1_H)$ in \cite [Proposition 2 (c)] {Ra94}, we have that
(iii) and (v) are equivalent. (iii) implies (iv) by \cite [Theorem
2.5.3] {Mo93}. Obviously, (iv) implies (v). $\Box$

\begin {Lemma} \label {4.3.2}  (cf.  \cite  {Mo93}, \cite  {Ra94} )
If $H$ is a finite dimensional Hopf algebra with
 $char  k > (dim { \ \ } H )^2 $,
then the following conditions are equivalent.

(i) $H$ is semisimple and cosemisimple.

(ii) $S_H^2 = id _H$.

(iii) $tr(S_H^2) \not=0.$
\end {Lemma}
{\bf Proof.}  Using formula $tr (S_H^2) = \epsilon _H (\Lambda _H^l)
\Lambda _{H^*} ^r (1_H)$ in \cite [Proposition 2 (c)] {Ra94}, we
have that (i) and (iii) are equivalent. (i) implies (ii) by \cite
[Theorem 2.5.3] {Mo93}. Obviously, (ii) implies (iii). $\Box$

\begin {Proposition}  \label {4.3.4}

Assume that  $D=A \stackrel {b} {\bowtie } H$ are finite dimensional
Hopf algebras.
  If $A  \stackrel {b} {\bowtie } H$ is   (co)semisimple,  then
$A$ and $H$ are  (co)semisimple.

\end {Proposition}

{\bf Proof .}  If  that $D$ is semisimple,  then $A$ and $H$ are
semisimple by Proposition \ref {4.3.3}. If  that $D$ is
cosemisimple, by Lemma \ref {4.3.4'}, $D^*$ is semisimple.
Considering $(A \stackrel {b} {\bowtie } H)^* \cong A^* \stackrel
{b} {\bowtie } H^*$ ( see \cite [Proposition 3.1.2]{Zh99}, note that
the evaluations in  \cite [Proposition 3.1.12]{Zh99} or \cite
[Proposition 1.11] {ZC01} are not the same in this paper ), we have
that $A^*$ and $H^*$ are semisimple by Proposition \ref {4.3.3}.
Consequently, $A$ and $H$ are cosemisimple.
 $\Box$

\begin {Theorem}  \label {4.3.5} Assume that
\  $\{ j_{A_1}, j_{A_2}, \cdots , j_{A_n} \}$ \ is a double
factorization of finite- dimensional Hopf algebra  $H$
 with respect to some method
$\sigma $ adding brackets. Then

(I)  $H$  is semisimple and cosemisimple  iff $A_i$ is semisimple
and cosemisimple for $i= 1, 2, \cdots , n, $  iff $tr (S_{A_i}^2)
\not=0$  for $i =1, 2 \cdots, n,$   iff $tr (S_{H}^2) \not=0$.

(II) If  $H$  is (co)semisimple, then  $A_i$ is (co)semisimple for
$i= 1, 2, \cdots , n;$

(III)
 If  $A_i$  is involutory for $i= 1, 2 ,
\cdots , n$,  then $H$ is  involutory.

(IV) $H$  admits a coquasitriangular structure iff  $A_i$ admits a
coquasitriangular
  structure
   for $i= 1, 2, \cdots , n.$

(V) If $char k =0$,
 then  the following are equivalent.

(i) $H$ is semisimple and cosemisimple; (i)$'$  $H$ is semisimple;
(i)$''$ $H$ is  cosemisimple.

(ii)  $A_i$  is semisimple and cosemisimple for $i=1, 2, \cdots, n$;
(ii)$'$ $A_i$  is semisimple for $i=1, 2, \cdots, n;$ (ii)$''$
$A_i$ is
 cosemisimple for $i=1, 2, \cdots, n.$

(iii)  $ A_i$  is involutory  for $i=1,2, \cdots n.$

(iv) $H$ is involutory.

(v) $tr (S_H ^2) \not= 0$.

(vi)  $tr (S_{ A_i}^2) \not= 0$  for $i=1,2, \cdots n.$

(VI) If  $char k > (dim H) ^2$, or $char k > (dim A_i)^2$ for $i=1,
2, \cdots n$, then  the following are equivalent.

(i) $H$ is semisimple and cosemisimple.

(ii)  $A_i$  is semisimple and cosemisimple for $i=1,2, \cdots ,n.$

(iii)  $ A_i$  is involutory  for $i=1,2, \cdots n.$

(iv) $H$ is involutory.

(v) $tr (S_H ^2) \not= 0$.

(vi)  $tr (S_{ A_i}^2) \not= 0$   for $i=1,2, \cdots n.$

\end{Theorem}
{\bf Proof.} By Theorem  \ref {4.2.2}, there exists $\{ \alpha _i ,
\beta _i \mid i = 1, 2, \cdots , n \}$   such that $D=: \sigma (A_1
{ \ } _{\alpha _1 } \bowtie _{\beta _1} {A_2} {\ }_{\alpha _2}
\bowtie _{\beta _2} {A_3}{ \ } _{\alpha _3} \bowtie _{\beta _3}
\cdots _{\alpha _{n-1}}\bowtie  _{\beta _{n-1}}A_n  )\cong H $
   \ \ \   (   as  Hopf algebras )   and
 the isomorphism is $m^{n-1}.$
 It is clear that
 \begin {eqnarray} \label {e31}S_D ^2 = S_{A_1} ^2 \otimes S_{A_2}^2 \otimes \cdots \otimes
 S_{A_n}^2 \end {eqnarray} (see \cite [Proposition 1.6]{ZC01}) and
\begin {eqnarray} \label {e32} tr (S_D ^2) = tr ( S_{A_1} ^2) tr(S_{A_2}^2) \cdots tr ( S_{A_n}^2) \end {eqnarray}
(see \cite [Theorem XIV.4.2]{Ka95}).

(I) If $H$ is semisimple and cosemisimple, then $tr(S_H ^2) \not=0,$
so $tr (S_{A_i}^2) \not=0$ by formula (\ref {e32}) for $i =1, 2
\cdots, n$. It follows from \cite [Proposition 2 (c)] {Ra94} that
$A_i$ is semisimple and cosemisimple for $i =1, 2, \cdots, n$.
Similarly, we can show the others.

(II) It follows from Proposition \ref {4.3.4}.

(III) It follows from  formula (\ref {e31}).

(IV) It follows from  the dual result of \cite [Corollary 2.3]
{Ch98} or \cite [Corollary $7.4.7^\circ $] {Zh99}.

(V) By Lemma \ref {4.3.1}, (i), (i)$'$ and (i)$''$ are equivalent;
(ii), (ii)$'$ and (ii)$''$  are equivalent; (ii), (iii) and (vi) are
equivalent; (i), (iv) and (v) are equivalent. By (I), (i),  (ii),
(v) and (vi) are equivalent.

(VI) By Lemma \ref {4.3.2}, (ii), (iii) and (vi) are equivalent. By
(I), (i), (ii), (v) and (vi) are equivalent.  By (III), (iii)
implies (iv).

 Now
we show that (iv) implies (v). Let $char  k = p $. Since $p$ is
prime and $dim  H = (dim  A_1)(dim  A_2)\cdots (dim  A_n)$ with $p
> ( dim  A_i )^2$, we have that $p $ does not divide $dim H$.
Therefore, $tr (S_H^2) = tr (id _H) \not=0$. $\Box$

{\bf Acknowledgement}\ \ The authors thank the referee and the
editor for valuable suggestions.

\end {document}